\documentclass[10pt,reqno,final]{article}

\usepackage{amsmath,amsfonts,amssymb,amsthm,version}
\usepackage{mathrsfs,fancybox,pifont}
\usepackage{graphicx}
\usepackage{url,hyperref}
\usepackage[notcite,notref]{showkeys}
\usepackage{color}
\usepackage{subfigure,multirow}
\usepackage{epstopdf}
\usepackage{cases}
\usepackage{mathtools}
\usepackage{algorithm,algorithmic}
\usepackage{authblk}
\usepackage{lipsum}

\allowdisplaybreaks

\theoremstyle{plain}

\theoremstyle{definition}

\theoremstyle{remark}

\title{A temporal multiscale method and its analysis for a system of fractional differential equations}
\author{Zhaoyang Wang\thanks{Department of Applied Mathematics, School of Mathematics and Physics, University of Science and Technology Beijing, Beijing 100083, China (zhaoyang584520@163.com)} \, and
Ping Lin\thanks{Corresponding author. Division of Mathematics, University of Dundee, Dundee DD1 4HN, United Kingdom (p.lin@dundee.ac.uk)}}

\affil{}

\date{\today}

\begin{document}
\maketitle

\begin{abstract}
In this paper, a nonlinear system of fractional ordinary differential equations with multiple scales in time is investigated. We are interested in the effective long-term computation of the solution. The main challenge is how to obtain the solution of the coupled problem at a lower computational cost. We analysize a multiscale method for the nonlinear system where the fast system has a periodic applied force and the slow equation contains fractional derivatives as a simplication of the atherosclerosis with a plaque growth. A local periodic equation is derived to approximate the original system and the error estimates are given. Then a finite difference method is designed to approximate the original and the approximate problems. We construct four examples, including three with exact solutions and one following the original problem setting, to test the accuracy and computational efficiency of the proposed method. It is observed that, the computational time is very much reduced and the multiscale method performs very well in comparison to fully resolved simulation for the case of small time scale separation. The larger the time scale separation is, the more effective the multiscale method is.


\medskip
\end{abstract}

\section{Introduction}
Many problems in nature have multiscale characteristics. For example, material weathering, material fracture due to atomistic defects such as impurities, chemical reaction of many substances mixed with fast and slow reaction rates, or some biological problems with a slow growth. In these processes large-scale or long-term phenomena are often influenced by important small-scale or short-term effects, and the variation of long-term properties will also effect the short-term response. A method that fully resolve short-term effects may not be effective for catching a long-term property. Therefore we must effectively deal with these two-way coupled problems. 

In the past two decades, people are interested in developing efficient multiscale methods for different temporal/spatial multiscale problems. For spatial multiscale problem, people mainly develop computable models such as atomistic-to-continuum coupling (AtC) models methods to simulate material behaviors involved with the atomistic scale, where the simulation of materials at the full atomistic scale is not possible due to its computatonal cost (See, e.g. \cite{Tadmor1996, Lin2003, Weinan2007, Shenoy1999, Weinan2006, Shapeev2011, Ortner2012, Abdulle2012, Luskin2013} and most recently \cite{Wang2018} for AtC of both simple and complex lattice structures). For temporal multiscale problem, one of the most prominent class of techniques is the Heterogeneous Multiscale Method (HMM) 
\cite{Engquist2005, E2003, Weinan2003, Abdulle2012(2)} that aims at efficient decoupling of a macro-scale and micro-scale problem, where the latter one enters the macro-scale problem in terms of temporal averages. In the context of ordinary differential equations, complexity plays a smaller role. The authors of \cite{Ariel2009, Ariel2009(2)} numerically realize analytical averaging techniques for problems that are also characterized by resonance effects. Frei and Richter addressed the Navier-Stokes flow problem with multiscale characteristics in time \cite{Frei2020, Richter2021(1)}. There are other applications and developments related to temporal multiscale problems (see, e.g. \cite{E2011, Richter2021(2), Trucu2013}). To our knowledge, none of those has considered fractional differential equations.

In this paper, we consider a temporal multiscale problem with fractional differentiations. The HMM \cite{Weinan2003, E2003, E2011, Abdulle2012, Engquist2005} is a general framework for designing multiscale algorithms, and in particular, temporal multiscale problems. Here, we consider a nonlinear coupled system of fractional ODEs, which is simplified from describing the problem of atherosclerosis \cite{Frei2020, Yang2016}. In \cite{Frei2020}, a multiscale approach has been rigirously developed for a linear form of ODEs as a simplification of the evolution of plaques in blood vessels, where a periodic force is applied to the flow due to the periodic heart pulse nature. Considering that the plaque accumulation or evolution may have memory effects in time and that a blood flow may be highly nonlinear, we thus aim to develop the method for a fractional order multiscale nonlinear ODE system with a periodic applied force. This development will biuld up a fundation to deal with many types of fractional order multiscale ODEs/PDEs in the future. 

The study of fractional calculus has attracted the attention of many scholars in the last 20 years, and it has been used in the modeling of problems in science, such as fluid mechanics \cite{Song2016}, biomathematics \cite{Kosztolowicz2020}, anomalous diffusion \cite{Zhaoyang2021, Sandev2019, Kilbas2006}, etc. On the contrary to the integer one, the fractional derivative has memory effects and non-local characteristics, which naturally appears in various problems in real world and leads to an increase in computational cost \cite{Gao2011}. In many cases, we are interested in the behavior of solution in a long period of if time. For a temporal multiscale problem with macroscopic and microscopic features in time, we simply increase the time step size for a long time computation, the accuracy can not be guaranteed and the temporal microscopic behaviour is not visible under a large time step. Our basic idea is to construct a local periodic problem to approximate the microscopic bahaviour of the original problem. The advantage is that the fast variable with microscopic feature in each macro-step can be solved effectively and independently. We prove that the constructed local periodic problem can effectively approximate the original problem by analyzing the error of the coupled micro-macro multiscale method. Numerical tests are implemented to show the correctness of our algorithm. By using numerical schemes proposed, we find a significant acceleration effect compared with the fully resolved simulation at least when the time scale separation is larger.

The rest of this manuscript is organized as follows. In Section \ref{section2}, we introduce a temporal multiscale problem with fractional derivatives considered and prove several auxiliary lemmas to obtain effective approximate periodic equations. In section \ref{section3}, we conduct the error analysis at the continuous level. In Section \ref{section4}, we present the temporal discretisation of the multiscale scheme for solving the proposed problem. Numerical tests are presented in Section \ref{section5} to illustrate the accuracy and computational efficiency of our method. Some conclusions and remarks are provided in Section \ref{section6}.

\section{Model problem and preliminaries}
\label{section2}

We consider a nonlinear coupled system of fractional ordinary differential equations.

\ \noindent\textit{Problem 1}. \textit{On I=[0.T], let $v$, $u$$\in C^1[0,T]$,
	\begin{equation}\label{the1}
		\begin{split}
			&\frac{d}{dt}v(t)+g(u(t),v(t))=f(t)\\
			&\frac{d^\alpha}{dt^\alpha}u(t)=\varepsilon R(t,u(t),v(t))\\
			&v(0)=v_0,\quad u(0)=u_0,\quad t\in I.
		\end{split}
\end{equation}}
Here, $v(t)$ and $u(t)$ represent the fast variable and the slow variable, respectively. $\frac{\partial^\alpha}{\partial t^\alpha}$ is the Caputo fractional derivative of order $0<\alpha<1$ \cite{Podlubny1999}. $\varepsilon\ll 1$ is a small parameter that controls the change of $u$. The force term $f(t)$ is a local periodic function with period 1 without loss of generality, that is, $f([t])=f([t]+1)$, $[t]$ is the integer ceiling function.

In addition to describing the chemical kinetics system \cite{Singh2017}, when $g=\lambda(u(t))v(t)$, problem 1 is a simplified atherosclerosis model with a certain type of memory effect included in the vessel wall or plaque growth. Macrophages in the vessel wall take up low density lipoproteins (LDL), which carry cholesterol and triglycerides to the tissues, and are finally transformed into foam cells, which are engorged with lipids \cite{Hahn2009}. The flow dynamics in the vessels are on the scale of milliseconds to seconds. However, the growth of plaque is usually described in terms of months, which is a slow process. For the formulation of the plaque growth equation, we consider the existence of a variety of cells in the vessel wall such as endothelial cells and vascular smooth muscle cells \cite{Yang2016,Libby2002}, and thus the diffusion of macrophages and foam cells is most likely to be anomalous. It may thus be more suitable to use a fractional operator representation.

\ \noindent\textit{Remark 1.1}. Problem \ref{the1} is a simplification of the problem of atherosclerosis described by a coupled PDEs. The goal of this paper is to present an idea to deal with the atherosclerosis problem based on this nonlinear fractional ODE model. The idea may be applied to the PDE problem, though analysis and implementation would be much more complicated and considered in future (Preliminary results have been obtained).\\

To apply the multiscale method and ensure the existence of the solution, it is necessary to have some assumptions about $R(t,u,v)$ and $g(u,v)$.

\ \noindent\textit{Assumption 2}. \textit{Let $v$, $u$$\in C^1[0,T]$. $R(t,u,v)$ and its derivative are bounded
	\begin{equation}\label{the2}
	\begin{split}
	\max \limits_{t\in[0,T]}\left\{|R(t,u,v)|,\left|\frac{\partial R(t,u,v)}{\partial t}\right|\right\}\leq C_R,
	\end{split}
	\end{equation}
	and Lipschitz continuous with respect to slow and fast variables
	\begin{equation}\label{the3}
	\begin{split}
	|R(t,u_1,v_1)-R(t,u_2,v_2)|\leq L_1|u_1-u_2|+L_2|v_1-v_2|. 
	\end{split}
	\end{equation}}

\ \noindent\textit{Assumption 3}. \textit{Let $v$, $u$$\in C^1[0,T]$. $g(u,v)$ is bounded
	\begin{equation}\label{the4}
	\begin{split}
	|g(u,v)|\leq C_g,
	\end{split}
	\end{equation}
	and $g(u,v)$ is differentiable with a bounded derivative
	\begin{equation}\label{the5}
	\begin{split}
	\left|\frac{\partial g(u,v)}{\partial u}\right|\leq L_3, \ 0<g_{min}\leq\frac{\partial g(u,v)}{\partial v}\leq L_4.
	\end{split}
	\end{equation}}

\ \noindent\textit{Remark 3.1}. From Eq. (\ref{the5}), we can obtain the Lipschitz condition with respect to slow and fast variables
\begin{equation}\label{the6}
\begin{split}
|g(u_1,v_1)-g(u_2,v_2)|\leq L_3|u_1-u_2|+L_4|v_1-v_2|. 
\end{split}
\end{equation}

\ \noindent\textit{Remark 3.2}. For the linear form given in \cite{Frei2020}, i.e., $g(u,v)=\lambda(u)v$, assumption 3 is reduced to $\left|\frac{d\lambda(u)}{du}\right|\leq C_\lambda$. (\ref{the5}) ensures that problem \ref{the1} has a local periodicity. Furthermore, we would like to point out that assumption 3 is reasonable if this model is used to describe the simplified atherosclerosis problem \cite{Frei2020}. Its derivative is bounded as long as the growth of plaque does not change the topology of the computing domain.\\

We introduce a new variable
\begin{equation}\label{the7}
\begin{split}
U(t)=u_0+I^{\alpha}_{0^+}\int_{s}^{s+1}\frac{d^\alpha u(r)}{dr^\alpha}dr,
\end{split}
\end{equation}
where $I^{\alpha}_{0^+}$ is the Riemann-Liouville fractional integral operator and defined by
\begin{equation}\label{the8}
\begin{split}
I^{\alpha}_{0^+}u(t)=\frac{1}{\Gamma(\alpha)}\int_{0}^{t}(t-s)^{\alpha-1}u(s)ds.
\end{split}
\end{equation}
By inserting $R(t,U(t),v(s))$, we have
\begin{equation}\label{the9}
\begin{split}
&\frac{d^{\alpha}U(t)}{dt^{\alpha}}=\int_{t}^{t+1}\frac{d^\alpha u(s)}{ds^\alpha}ds=\varepsilon\int_{t}^{t+1}R(t,u(s),v(s))ds\\
&=\varepsilon\int_{t}^{t+1}R(t,U(t),v(s))ds+\varepsilon\int_{t}^{t+1}\left(R(t,u(s),v(s))-R(t,U(t),v(s))\right)ds. 
\end{split}
\end{equation}

\ \noindent LEMMA 4. \textit{Let $u\in C^1[0,T]$, $v\in C^1[0,T]$, and let assumption 2 be satisfied. It holds that
	\begin{equation}\label{the10}
	\begin{split}
	\left|\varepsilon\int_{t}^{t+1}\left(R(t,u(s),v(s))-R(t,U(t),v(s))\right)ds\right|\leq C \varepsilon,
	\end{split}
	\end{equation}
	with a constant $C=4L_1|u_0|$.}\\
\begin{proof} 
	Applying the fractional integral operator $I^\alpha$ on both sides of the second equation of (\ref{the1}), we obtain
	\begin{equation}\label{the11}
	\begin{split}
	|u(t)|&=\varepsilon\left|u_0+\int_{0}^{t}\frac{(t-s)^{\alpha-1}}{\Gamma(\alpha)}R(t,u(s),v(s))ds\right|\\
	&\leq|u_0|\varepsilon+\frac{C_RT^{\alpha}}{\Gamma(\alpha+1)}\varepsilon.
	\end{split}
	\end{equation}
	Similarly, we can also get an estimate of the new variable $U(t)$
	\begin{equation}\label{the12}
	\begin{split}
	|U(t)|\leq |u_0|+\frac{\varepsilon C_RT^{\alpha}}{\Gamma(\alpha+1)}.
	\end{split}
	\end{equation}
	By Lipschitz continuity of $R$, we have 
	\begin{equation}\label{the13}
	\begin{split}
	&\left|\varepsilon\int_{t}^{t+1}\left(R(u(s),v(s))-R(U(t),v(s))\right)ds\right|\leq\varepsilon L_1\int_{t}^{t+1}\left|U(t)-u(s)\right|ds\\
	&\leq \varepsilon L_1|U(t)|+\varepsilon L_1\int_{t}^{t+1}|u(s)|ds\leq C\varepsilon.
	\end{split}
	\end{equation}
\end{proof}

We can thus approximate the averaged evolution equation for $D^{\alpha}U$ by
\begin{equation}\label{the14}
\begin{split}
\frac{d^{\alpha}U(t)}{dt^{\alpha}}=\varepsilon\int_{t}^{t+1}R(t,U(t),v(s))ds+O(\varepsilon).
\end{split}
\end{equation}
The discretization of Eq. (\ref{the14}) in a macro-time step $T_n\to T_{n+1}=T_n+\triangle T$ is not effective because it is involved with the dynamic evolution of $v(T_n)$ to $v(T_{n+1})$ on the fast scale. We need deal with and approximate the fast variable solution $v(t)$ in a more effective way. We thus construct an auxiliary problem on a period to approximate the original problem \ref{the1}. Next, we give two lemmas to derive an equivalent local periodic problem, which may provide a fast scale solution $v_{U_n}$ without initial value. \\

\ \noindent LEMMA 5. \textit{Let $v\in C^1[0,T]$, $\lambda\in C^{1}[0,T]$ and $f\in C[0,T]$. For the first equation of (\ref{the1}), we have the following estimate
	\begin{equation}\label{the15}
	\begin{split}
	\frac{d^{\alpha}U(t)}{dt^{\alpha}}=\varepsilon\int_{t}^{t+1}R(t,U(t),v_{U(t)}(s))ds+O(\varepsilon).
	\end{split}
	\end{equation}}
\begin{proof} 
	The mild solution of the first equation of (\ref{the1}) is
	\begin{equation}\label{the16}
	\begin{split}
	v(t)=v_0+\int_{0}^{t}f(s)ds-\int_{0}^{t}g(u(s),v(s))ds.
	\end{split}
	\end{equation}
	Thus, we have the following estimate
	\begin{equation}\label{the17}
	\begin{split}
	|v(t)|\leq v_0+\left(\Vert f\Vert_{L^{\infty}[0,T]}+C_g\right)T.
	\end{split}
	\end{equation}
	We approximate the averaged fractional equation (\ref{the14}) by inserting $R(U(t),v_{U(t)}(s))$ for a fixed $U(t)$
	\begin{equation}\label{the18}
	\begin{split}
	&\frac{d^{\alpha}U(t)}{dt^{\alpha}}=\varepsilon\int_{t}^{t+1}R(t,U(t),v_{U(t)}(s))ds\\
	&+\varepsilon\int_{t}^{t+1} \left(R(t,U(t),v(s))-R(t,U(t),v_{U(t)}(s))\right)ds+O(\varepsilon).
	\end{split}
	\end{equation}
	For the second remainder, we have
	\begin{equation}\label{the19}
	\begin{split}
	&\varepsilon\int_{t}^{t+1} \left(R(t,U(t),v(s))-R(t,U(t),v_{U(t)}(s))\right)ds\\
	&\leq \varepsilon L_2 \int_t^{t+1}|v(s)-v_{U(t)}(s)|ds\leq 2L_2\left[v_0+\left(\Vert f\Vert_{L^{\infty}[0,T]}+C_g\right)T\right]\varepsilon.
	\end{split}
	\end{equation}
	The proof is completed.
\end{proof}

\ \noindent LEMMA 6. \textit{Let $u\in C^1[0,T]$, $v\in C^1[0,T]$ and $f([t])=f([t]+1)$. For a fixed $U$, if $v_U(0)=v_U(1)$ or $t\to\infty$, the periodic condition to the equation given below holds
	\begin{equation}\label{the20}
	\begin{split}
	\frac{d}{dt}v_{U}(t)+g(U,v_{U}(t))=f(t), \ v_U([t])=v_U([t]+1).
	\end{split}
	\end{equation}}
\begin{proof} 
	Let $w(t)=v_U([t]+1)-v_U([t])$, we have 
	\begin{equation}\label{the21}
	\begin{split}
	\frac{d}{dt}w(t)+g(U,v_U([t]+1))-g(U,v_U([t]))=0,
	\end{split}
	\end{equation}
	By the mean value theorem, we obtain the following
	\begin{equation}\label{the22}
	\begin{split}
	\frac{d}{dt}w(t)+\frac{dg(U,\xi)}{dv_{U}}w(t)=0, \ \xi\in(v_U([t],v_U([t]+1)),
	\end{split}
	\end{equation}
	with the solution
	\begin{equation}\label{the23}
	\begin{split}
	|w(t)|=|w(0)|e^{-\frac{dg}{dv_{U}}t}\leq |w(0)|e^{-g_{min}t}.
	\end{split}
	\end{equation}
	Hence, we obtain $|w(t)|=0$. 
\end{proof} 

Based on the above lemma, for this temporal multiscale problem, we construct the following auxiliary problem.\\

\ \noindent \textit{Problem 7.} \textit{On I=[0.T], let $v$, $u$$\in C^1[0,T]$,
	\begin{equation}\label{the24}
	\begin{split}
	&\frac{d}{dt}v_U(t)+g(U,v_U(t))=f(t)\\
	&\frac{d^\alpha}{dt^\alpha}U(t)=\varepsilon\int_t^{t+1}R(t,U(t),v_{U(t)}(s))ds\\
	&v_U([t])=v_U([t]+1),\quad U(0)=u_0.
	\end{split}
	\end{equation}}

\ \noindent\textit{Remark 7.1}. If the assumption of $f(t)$ is periodic, i.e., $f(t)=f(t+1)$, the periodic condition of Eq. (\ref{the24}) should be $v_U(0)=v_U(1)$ \cite{Frei2020}. 

\ \noindent\textit{Remark 7.2}. Lemma 6 implies that $t\to\infty$, $v(t)$ will show local periodic properties. Since we focus on long-term simulation, the initial non-periodic part has little effect on it, which is confirmed in numerical experiments.\\

\section{Error analysis of the fractional multiscale problem}
\label{section3}
In this section, we will prove that the difference between variable $U(t)$ in the auxiliary problem \ref{the24} and variable $u(t)$ in original problem \ref{the1} is small. The error estimate between $v_{U(t)}(t)$ and $v(t)$ is also obtained. Therefore, the solution of problem \ref{the24} can be used as an effective approximation of the solution of problem \ref{the1}.\\

\ \noindent LEMMA 8. \textit{Let $u$ be fixed and $f\in L^\infty[0,T]$. For the local periodic problem
	\begin{equation}\label{the25}
	\begin{split}
	\frac{d}{dt}v_{u}(t)+g(u,v_u(t))=f(t), \ v_u([t])=v_u([t]+1).
	\end{split}
	\end{equation}
	It holds 
	\begin{equation}\label{the26}
	\begin{split}
	\left|\frac{\partial v_u(t)}{\partial u}\right|\leq C_{L8}.
	\end{split}
	\end{equation}
	where $C=\{g_{min}, L_3, L_4\}$.}\\
\begin{proof} 
	Differentiate both sides of Eq. (\ref{the25}) with respect to $u$. Let $p(t)=\frac{\partial v_u(t)}{\partial u}$, we have
	\begin{equation}\label{the27}
	\begin{split}
	\frac{\partial}{\partial t}p(t)+\frac{\partial g(u,v_u(t))}{\partial u}+\frac{\partial g(u,v_u(t))}{\partial v_u}p(t)=0.
	\end{split}
	\end{equation}
	Since $u$ is a parameter, the boundary value condition is local periodic
	\begin{equation}\label{the28}
	\begin{split}
	p([t])=p([t]+1).
	\end{split}
	\end{equation}
	Multiplying $exp\{\int_{[t]}^t\frac{\partial g(u,v_u(s))}{\partial v_u}ds\}$ on both sides for Eq. (\ref{the27}) and integrating in $t$, we obtain
	\begin{equation}\label{the29}
	\begin{split}
	p(t)=p([t])e^{-\int_{[t]}^t\frac{\partial g(u,v_u(s))}{\partial v_u}ds}-e^{-\int_{[t]}^t\frac{\partial g(u,v_u(s))}{\partial v_u}ds}\int_{[t]}^t\frac{\partial g(u,v_u(s))}{\partial u}e^{\int_{[t]}^s\frac{\partial g(u,v_u(r))}{\partial v_u}dr}ds.
	\end{split}
	\end{equation}
	Using the local periodic boundary conditon, we have
	\begin{equation}\label{the30}
	\begin{split}
	|p([t])|&=\frac{e^{-\int_{[t]}^{[t]+1}\frac{\partial g(u,v_u(s))}{\partial v_u}ds}}{1-e^{-\int_{[t]}^{[t]+1}\frac{\partial g(u,v_u(s))}{\partial v_u}ds}}\left|\int_{[t]}^{[t]+1}\frac{\partial g(u,v_u(s))}{\partial u}e^{\int_{[t]}^s\frac{\partial g(u,v_u(r))}{\partial v_u}dr}ds\right|\\
	&\leq \frac{L_3e^{L_4}}{1-e^{-g_{min}}}.
	\end{split}
	\end{equation}
	For $\forall t\in[[t],[t]+1]$
	\begin{equation}\label{the31}
	\begin{split}
	|p(t)|&\leq \frac{L_3e^{L_4}}{1-e^{-g_{min}}}+\int_{[t]}^t\left|\frac{\partial g(u,v_u(s))}{\partial u}e^{\int_{[t]}^s\frac{\partial g(u,v_u(r))}{\partial v_u}dr}\right|ds\\
	&\leq \frac{L_3e^{L_4}}{1-e^{-g_{min}}}+L_3e^{L_4}:=L_8.
	\end{split}
	\end{equation}
\end{proof}

The following lemma will give the error estimate between the local periodic solution $v_{u(t)}(t)$ and the fast variable $v(t)$. For a fixed $u(t)$, we consider the family of periodic solutions
\begin{equation}\label{the32}
\begin{split}
&\frac{\partial}{\partial s}v_{u(t)}(s)+g(u(t),v_{u(t)}(s))=f(s).\\ 
&\ v_{u(t)}([s])=v_{u(t)}([s]+1). \ s\in[[s],[s]+1] \ and \ t\in[0,T].
\end{split}
\end{equation}

\ \noindent LEMMA 9. \textit{Let $u, v\in C^1[0,T]$, assumptions 2 and 3 hold. For the initial values to (\ref{the1}) and (\ref{the32}) agree, i.e., $v_{u(0)}(0)=v_0$. It holds
	\begin{equation}\label{the33}
	\begin{split}
	|v_{u(t)}(t)-v(t)|\leq C_{L9}\varepsilon,
	\end{split}
	\end{equation}
	with a constant $C=\{C_R,C_{L8},L_4,\alpha\}$.}\\
\begin{proof}
	For $v_{u(t)}(t)$, we have
	\begin{equation}\label{the34}
	\begin{split}
	\frac{\partial}{\partial t}v_{u(t)}(t)=\frac{\partial v_{u(t)}(t)}{\partial s}+\frac{d v_{u(t)}}{d u(t)}u'(t),
	\end{split}
	\end{equation}
	the function $v_{u(t)}(t)$ fulfils the following equation
	\begin{equation}\label{the35}
	\begin{split}
	\frac{\partial}{\partial t}v_{u(t)}(t)-\frac{d v_{u(t)}}{d u(t)}u'(t)+g(u(t),v_{u(t)}(t))=f(t), \ v_{u(0)}(0)=v_0.
	\end{split}
	\end{equation}
	Let $w(t)=v_{u(t)}(t)-v(t)$, it holds
	\begin{equation}\label{the36}
	\begin{split}
	\frac{\partial}{\partial t}w(t)=g(u(t),v(t))-g(u(t),v_{u(t)}(t))+\frac{d v_{u(t)}}{d u(t)}u'(t), \ w(0)=0.
	\end{split}
	\end{equation}
	For $u'(t)$, we use the fractional differential operator $^{C}D^{1-\alpha}$ on both sides of the second equation of (\ref{the1}), we have 
	\begin{equation}\label{the37}
	\begin{split}
	\left|\frac{d}{dt}u(t)\right|=\frac{\varepsilon}{\Gamma(\alpha)}\left|\int_0^{t}(t-s)^{\alpha-1}\frac{dR(s,u,v)}{ds}ds\right|\leq \frac{C_{R}T^{\alpha}\varepsilon}{\Gamma(\alpha+1)}.
	\end{split}
	\end{equation}
	To estimate the right term in (\ref{the36}), we integrate (\ref{the36}) in $t$ and use the Lipschitz condition (\ref{the6}) and lemma 8, we obtain
	\begin{equation}\label{the38}
	\begin{split}
	|w(t)|&\leq w(0)+\int_0^{t}|g(u(t),v_{u(t)}(s))-g(u(t),v(s))|ds+C_{L8}\int_0^{t}|u'(s)|ds\\
	&\leq L_4\int_0^t|w(s)|ds+\frac{C_{L8}C_{R}T^{\alpha+1}\varepsilon}{\Gamma(\alpha+1)}
	\end{split}
	\end{equation}
	By Gronwall inequality, we can easily get the following estimates
	\begin{equation}\label{the39}
	\begin{split}
	|w(t)|\leq \frac{C_{L8}C_{R}T^{\alpha+1}}{\Gamma(\alpha+1)}e^{L_{4}T}\varepsilon:=C_{L9}\varepsilon.
	\end{split}
	\end{equation}
\end{proof} 

Lemma 9 proves that the initial value problem \ref{the1} can be transformed into a local periodic problem \ref{the24} with fixed slow variable $u(t)$. Next we will prove the last lemma to show the final error estimate between the solution of the original problem \ref{the1} and the solution of the multiscale approximation problem \ref{the24}.\\

\ \noindent LEMMA 10. \textit{Let $u, v\in C^1[0,T]$, and assumptions 2 and 3 hold. $(u(t),v(t))$ be defined by (\ref{the1}) and $(U(t),v_{U(t)})$ be defined by (\ref{the24}), with the initial condition $v_{U(0)}(0)=v_0$. For $t\in[0,T]$ and $T=O(\varepsilon^{-\alpha^{-1}})$, it holds 
	\begin{equation}\label{the40}
	\begin{split}
	|U(t)-u(t)|\leq C_{L10a}\varepsilon,\quad |v_{U(t)}(t)-v(t)|\leq C_{L10b}\varepsilon.
	\end{split}
	\end{equation}}
\begin{proof} 
	Let $w(t)=U(t)-u(t)$, we have 
	\begin{equation}\label{the41}
	\begin{split}
	\frac{d^\alpha}{dt^\alpha}w(t)&=\varepsilon\int_t^{t+1}R(t,U(t),v_{U(t)}(s))ds-\varepsilon R(t,u(t),v(t))\\
	&=\varepsilon\int_t^{t+1}\left(R(t,U(t),v_{U(t)}(s))-R(t,u(t),v(t))\right)ds
	\end{split}
	\end{equation}
	Lipschitz continuity of $R$ in assumption 2 gives
	\begin{equation}\label{the42}
	\begin{split}
	&\left|\frac{d^\alpha}{dt^\alpha}w(t)\right|\leq \varepsilon L_1\int_t^{t+1}|U(t)-u(t)|ds+\varepsilon L_2\int_t^{t+1}|v_{U(t)}(s)-v(t)|ds\\
	&\leq \varepsilon L_1|w(t)|+\varepsilon L_2\int_t^{t+1}\left(|v_{U(t)}(s)-v_{u(t)}(s)|+|v_{u(t)}(s)|+|v(t)|\right)ds.
	\end{split}
	\end{equation}
	By using lemma 4 and lemma 8, we get the following estimates
	\begin{equation}\label{the43}
	\begin{split}
	&\int_t^{t+1}|v_{U(t)}(s)-v_{u(t)}(s)|+|v_{u(t)}(s)|+|v(t)|ds\\
	&\leq \int_t^{t+1}(C_{L8}|w(t)|+2C)ds\\
	&=C_{L8b}|w(t)|+2C.
	\end{split}
	\end{equation}
	We combine (\ref{the42}) and (\ref{the43}) to obtain
	\begin{equation}\label{the44}
	\begin{split}
	&-C\varepsilon w(t)-C\varepsilon\leq \frac{d^\alpha}{dt^\alpha}w(t)\leq C\varepsilon w(t)+C\varepsilon\\
	&\Rightarrow |w(t)|\leq \frac{C\varepsilon}{\Gamma(\alpha)}\int_0^t (t-s)^{\alpha-1}|w(s)|ds+C\varepsilon,
	\end{split}
	\end{equation}
	where $C=\{C_{L8},C_{g},v_0,\Vert f\Vert_{L^\infty},\alpha\}$. Applying Gronwall's inequality technique it is easy to get
	\begin{equation}\label{the45}
	\begin{split}
	|w(t)|\leq C\varepsilon e^{C\varepsilon t^{\alpha}},
	\end{split}
	\end{equation}
	which satisfies $|w(t)|\leq C_{L10a}\varepsilon$ for $t\leq T=O(\varepsilon^{-\alpha^{-1}})$. Finally, we insert $v_{u(t)}(t)$ and use lemma 9 to obtain
	\begin{equation}\label{the46}
	\begin{split}
	|v_{U(t)}(t)-v(t)|&\leq |v_{U(t)}(t)-v_{u(t)}(t)|+|v_{u(t)}(t)-v(t)|\\
	&\leq C|w(t)|+C_{L9}\varepsilon\leq C_{L10b}\varepsilon.
	\end{split}
	\end{equation}
\end{proof} 

We have completed the theoretical analysis of an approximate multiscale problem, which is the fundation for the numerical calculation of fractional temporal multiscale problems.

\section{Numerical algorithms}
\label{section4}
In this section, we give numerical algorithms for solving problems \ref{the1} and \ref{the24}. The numerical discretization is based on the finite difference method. Define $t_i=i\triangle t \ (i=0,1,2,...,N)$, where $\triangle t=T/N$ is the micro-scale time step.

For Eq. (\ref{the1}), the L1 approximation \cite{Gao2012} for the Caputo fractional derivative of order $0<\alpha<1$ is given by
\begin{equation}\label{the47}
\begin{split}
\frac{d^\alpha}{dt^\alpha}u(t_i)=\frac{(\triangle t)^{-\alpha}}{\Gamma(2-\alpha)}\left[a_0u_i-\sum\limits_{j=1}^{i-1}(a_{i-j-1}-a_{i-j})u_j-a_{i-1}u_0\right]+O(\triangle t^{2-\alpha}),
\end{split}
\end{equation}
where $a_0=1$, $a_j=(j+1)^{1-\alpha}-j^{1-\alpha}$ $(j\geq1)$.

Then, the governing equation (\ref{the1}) is discretized as follows
\begin{equation}\label{the48}
\begin{split}
&\frac{v_i-v_{i-1}}{\triangle t}+g(u_i,v_i)=f(t_i),\\
&\frac{(\triangle t)^{-\alpha}}{\Gamma(2-\alpha)}\left[a_0u_i-\sum\limits_{j=1}^{i-1}(a_{i-j-1}-a_{i-j})u_j-a_{i-1}u_0\right]=\varepsilon R(t_{i-1},u_{i-1},v_{i-1}).
\end{split}
\end{equation}
Here, for the calculation of the fast variable $v(t)$, we will use the first-order explicit/implicit Euler scheme for discretization. The calculation process is shown as follows

\begin{algorithm}
	\caption{Direct solution scheme.}
	\label{algorithm11}
	\begin{algorithmic}
		\STATE{Let $u(0)=u_0$, $v(0)=v_0$. Iterate for $i=1,2,...,N$}
		\FOR{$i=1:N$}
		\STATE{Calculate $u_i$ through Eq. $(\ref{the48})_2$}
		\STATE{Calculate $v_i$ through Eq. $(\ref{the48})_1$}
		\ENDFOR
	\end{algorithmic}
\end{algorithm}

For the calculation of Eq. (\ref{the24}), we first consider the identification of periodic solutions in a period. For a given initial value $v_{U}(0)$, iterate to $v_{U}(1)$, and then update the initial value to $v_{U}(1)$ and enter the loop until the given error limit is reached. This scheme is more general than the averaging acceleration scheme applied in \cite{Richter2021(1)}. It is suitable for calculating fractional periodic problems but the calculation time is longer. Below we give the theoretical basis for this shooting method.

Let $v(t)$ and $v^p(t)$ be the trial solution and the periodic solution respectively. $v_0$ is the initial value for trial, for fixed $U$, it is not difficult to obtain the following estimates
\begin{equation}\label{the49}
\begin{split}
&|v(t)-v^p(t)|=|v_0-v^p(0)|e^{-g_{min}t}\leq |v_0-v^p(0)|\\
&\Rightarrow |v(1)-v^p(1)|=|v(1)-v^p(0)|\leq e^{-g_{min}}|v(0)-v^p(0)|.
\end{split}
\end{equation}
We can see that the longer the time, the better the identification efficiency of periodic solutions.

\begin{algorithm}
	\caption{Periodic solution identification scheme.}
	\label{algorithm12}
	\begin{algorithmic}
		\STATE{Given an initial value $v_{U}(0)$ as a trial and let $tol>0$ be a given tolerance. Iterate for $k=1,2,...K$.}
		\STATE{Step 1 Solve one cycle of $(\ref{the24})_1$ for $v_{U,k}$ with the initial $v_{U,0}=v_{U}(0)$.}
		\STATE{Step 2 Calculate the error between the initial value and the final value	
			$err:=|v_{U}(1)-v_{U}(0)|.$}
		\STATE{Step 3 If $err\leq tol$, stop.}
		\STATE{Step 4 Update the initial $v_{U}(0)=v_{U}(1)$, and go to Step 1.}
	\end{algorithmic}
\end{algorithm}

From lemma 6, we can see that the decay of the nonstationary solution to periodic solution may be slow, which depends on $g_{min}$. If necessary, we can ignore the case of nonstationary solutions and focus on the acceleration of periodic problems.

Based on Algorithm \ref{algorithm12}, we give the following multiscale method.

\begin{algorithm}[H]
	\caption{Fractional multiscale method.}
	\label{algorithm13}
	\begin{algorithmic}
		\STATE{We split the time interval $I=[0,T]$ into subintervals of equal size $T_0<T_1<\cdots<T_M$, $\triangle T=T_m-T_{m-1}=T/M$. Let $U_0=u_0$. Iterate for $m=1,2,...$.}
		\STATE{Step 1 For $U=U_{m-1}$ solve the local periodic problem \ref{the24} to obtain $v_{U,k}$ by Algorithm \ref{algorithm12}. $m=1,2,...M$.}
		\STATE{Step 2 Calculate the integral term in $(\ref{the24})_2$
			
			$R_{m-1}=\frac{1}{K+1}\sum\limits_{k=0}^{K}R(t_k,U_{m-1},v_{U,k})$.}
		\STATE{Step 3 Step forward $U_{m-1}\rightarrow U_{m}$ and go to step 1 with the L1 approximation
			$U_{m}=\Gamma(2-\alpha)(\triangle T)^{\alpha}\varepsilon R_{m-1}+\sum\limits_{j=1}^{m-1}(a_{m-j-1}-a_{m-j})U_{j}+a_{m-1}U_0$.}
	\end{algorithmic}
\end{algorithm}

\section{Numerical experiments}
\label{section5}
In this section, we carry out several numerical experiments to test the correctness and effectiveness of the proposed numerical algorithm. Specifically, we construct exact solutions to verify the correctness of fully resolved simulation and multiscale method, and then compare the calculation time and accuracy.
\subsection{Local test}
\label{section5.1}
We first consider conducting a numerical test in a time period to verify Algorithm \ref{algorithm11} and Algorithm \ref{algorithm12}, which ensures the multiscale method is locally correct.

\ \textbf{Example 1.} The coupled differential equations are as follow
\begin{equation}\label{the52}
\begin{split}
&\frac{d}{dt}v(t)+u(t)sin\left(v(t)\right)+2v(t)+1=f_1(t)\\
&\frac{d^\alpha}{dt^\alpha}u(t)=\varepsilon R_1(t,u(t),v(t))\\
&v(0)=0,\quad u(0)=\frac{1}{2},\quad t\in [0,6].
\end{split}
\end{equation}
The approximate local periodic equation is
\begin{equation}\label{the53}
\begin{split}
&\frac{d}{dt}v_U(t)+Usin(v_U(t))+2v_U(t)+1=f_1(t)\\
&v_U(0)=v_U(6), \ U=u_0=\frac{1}{2}, \quad t\in [0,6].
\end{split}
\end{equation}
Here, $f_1(t)=-2t^2+10t+7+\left(\varepsilon\frac{\Gamma(3-\alpha)}{2}t^2+\varepsilon\Gamma(2-\alpha)t+\frac{1}{2}\right)sin(-t^2+6t)$, $R_1(t,u,v)=t^{2-\alpha}+t^{1-\alpha}+uv-\left(-t^2+6t\right)\left(\varepsilon\frac{\Gamma(3-\alpha)}{2}t^2+\varepsilon\Gamma(2-\alpha)t+\frac{1}{2}\right)$.

This problem has the following exact solution
\begin{equation}\label{the54}
\begin{split}
&v(t)=-t^2+6t,\quad u(t)=\varepsilon\frac{\Gamma(3-\alpha)}{2}t^2+\varepsilon\Gamma(2-\alpha)t+\frac{1}{2}.
\end{split}
\end{equation}

\begin{figure}[H]
	\centering
	\includegraphics[scale=0.36]{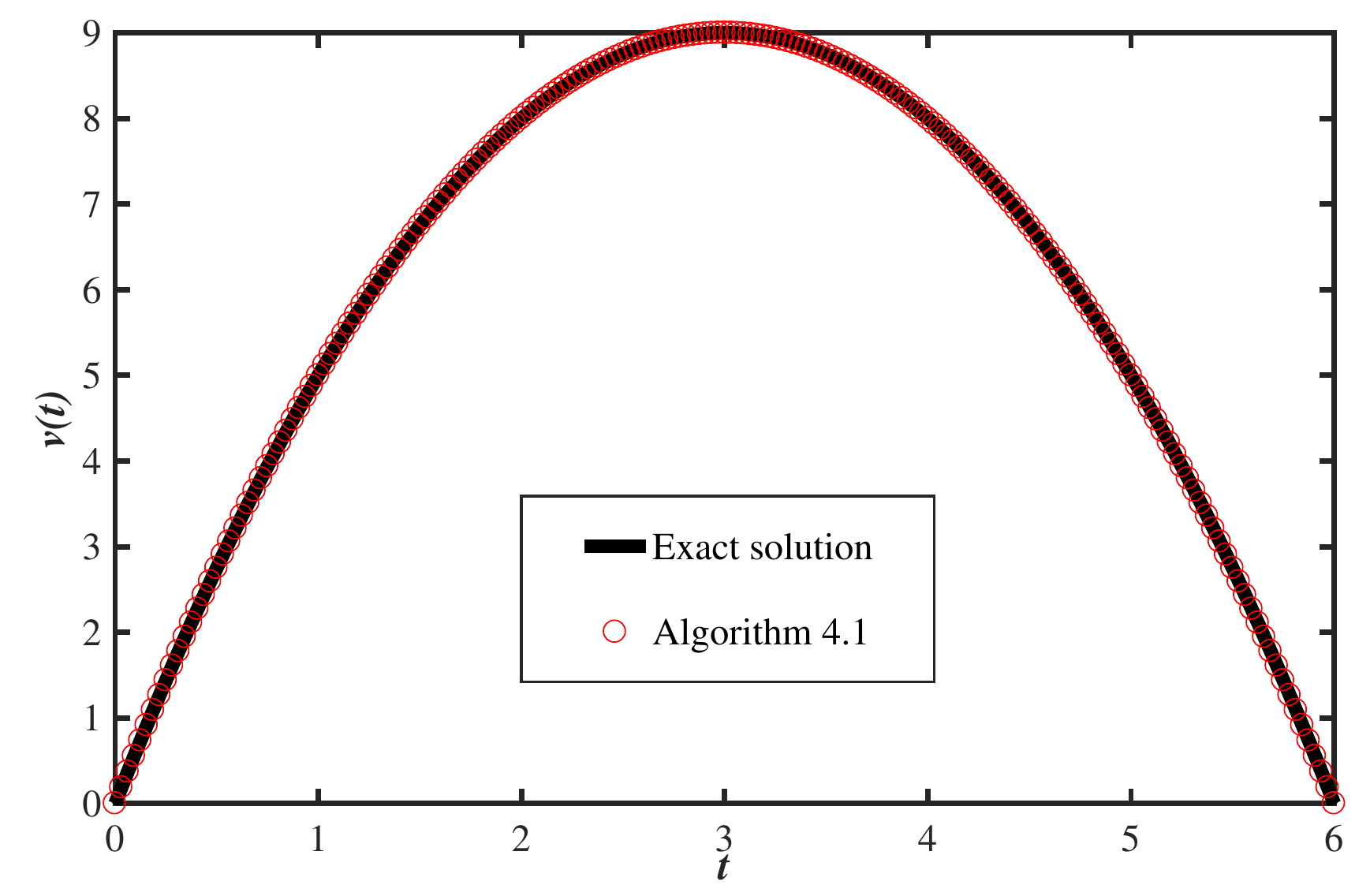}
	\includegraphics[scale=0.36]{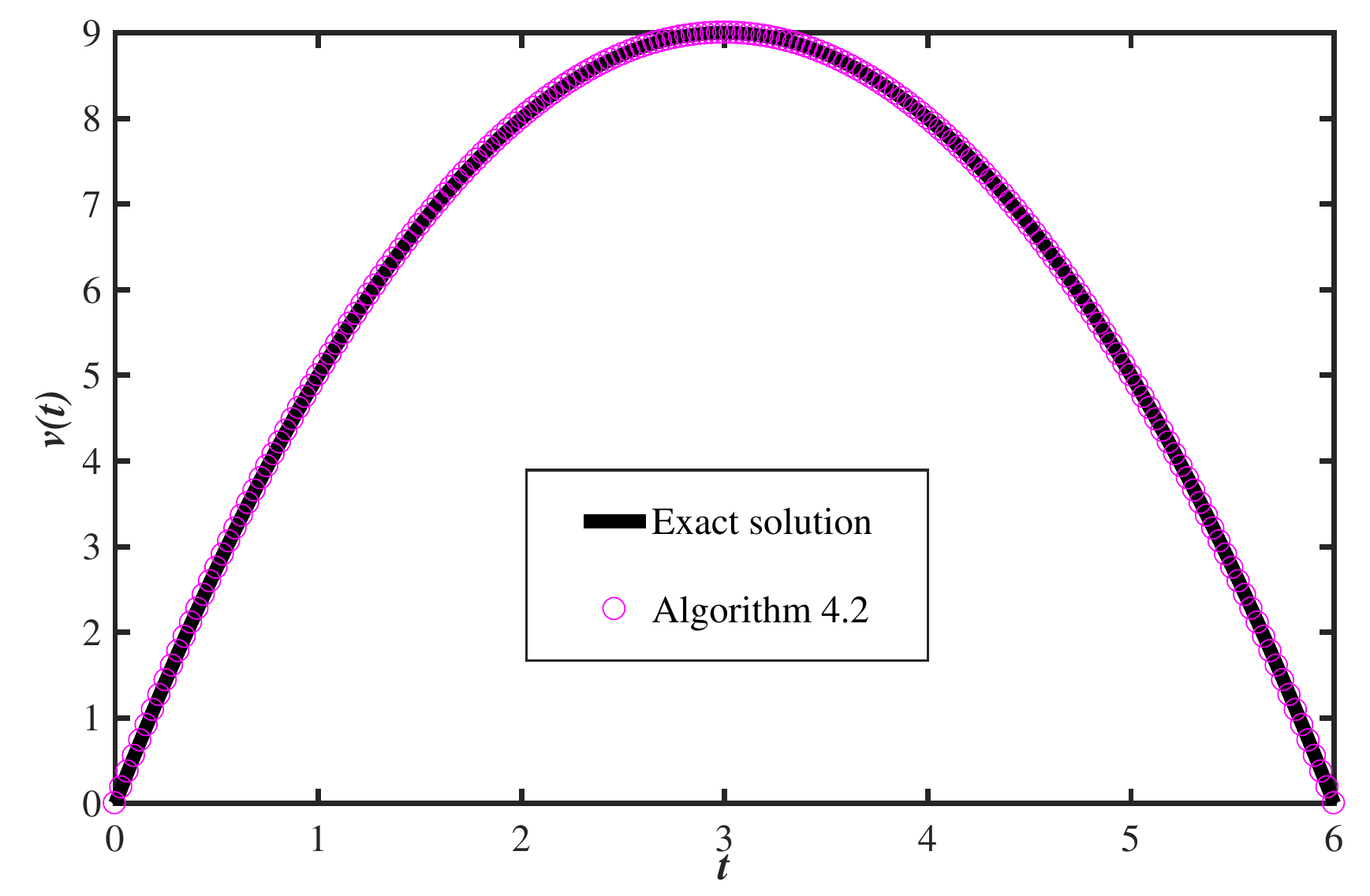}
	\caption{Comparison of the numerical computational results and exact solution.}
	\label{compare1}
\end{figure}

In the simulation, we set the time step $\triangle t=1/32$, $\alpha=0.6$ and $\varepsilon=5\times 10^{-5}$. The first-order explicit Euler scheme is used to discretize the fast variable equation. It can be seen from Fig. \ref{compare1} that the numerical computational results of Algorithm \ref{algorithm11} and Algorithm \ref{algorithm12} are in good agreement with the exact solution. In Algorithm \ref{algorithm12}, we take $tol=10^{-5}$. We would like to point out that Algorithm \ref{algorithm12} is also suitable for solving other types of boundary value problems (see, e.g., \cite{Ahmad2011, Ahmad2021}).

\subsection{Global test}
\label{section5.2}
We further consider a long term numerical simulation to test the effectiveness of the multiscale method. 

\ \textbf{Example 2.} The linear coupled system of fractional differential equations are given by
\begin{equation}\label{the55}
\begin{split}
&\frac{d}{dt}v(t)+\left(u(t)+1\right)v(t)=f_2(t)\\
&\frac{d^\alpha}{dt^\alpha}u(t)=\varepsilon R_2(t,u(t),v(t))\\
&v(0)=2,\quad u(0)=1,\quad t\in [0,10001],
\end{split}
\end{equation}
where $f_2(t)=sin(2\pi t)+2\pi tcos(2\pi t)+\left(\varepsilon\frac{\Gamma(2-\alpha)}{2}t^2+2\right)\left(tsin(2\pi t)+2\right)$, $R_2(t,u,v)=t^{2-\alpha}+\frac{\left(\varepsilon\frac{\Gamma(3-\alpha)}{2}t^2+1\right)\left(tsin(2\pi t)+2\right)}{uv}-1$.

The corresponding approximate equations are as follow
\begin{equation}\label{the56}
\begin{split}
&\frac{d}{dt}v_U(t)+(U+1)v_U(t)=f_2(t)\\
&\frac{d^\alpha}{dt^\alpha}U(t)=\varepsilon\int_t^{t+1}R_2(t,U(t),v_{U(t)}(s))ds\\
&v_U([t])=v_U([t]+1), \ U(0)=u_0=1,\quad t\in [0,10001].
\end{split}
\end{equation}
We has an exact solution for (\ref{the55}) 
\begin{equation}\label{the57}
\begin{split}
&v(t)=tsin(2\pi t)+2,\quad u(t)=\varepsilon\frac{\Gamma(3-\alpha)}{2}t^2+1.
\end{split}
\end{equation}

\begin{figure}[H]
	\centering
	\includegraphics[scale=0.35]{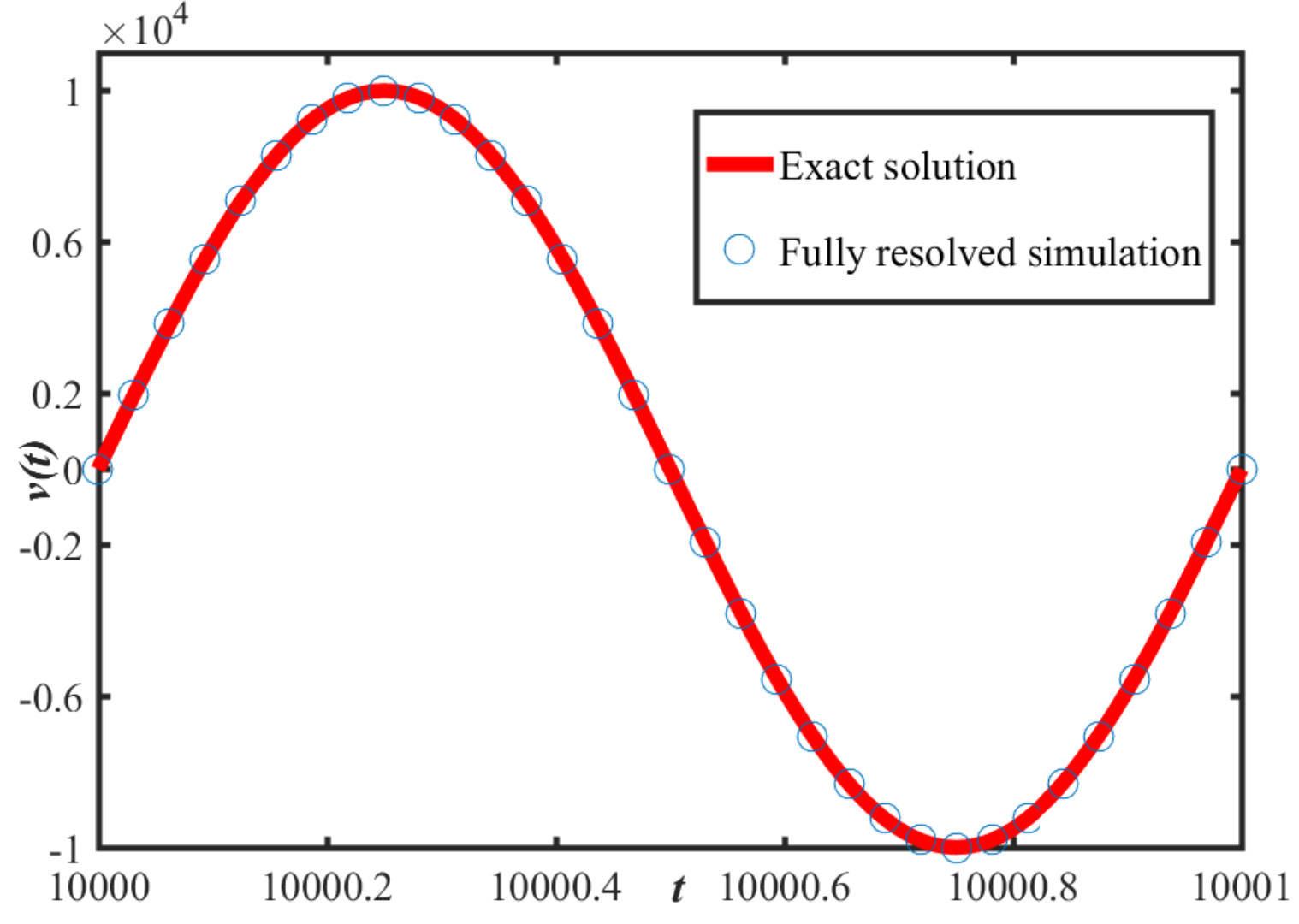}
	\includegraphics[scale=0.35]{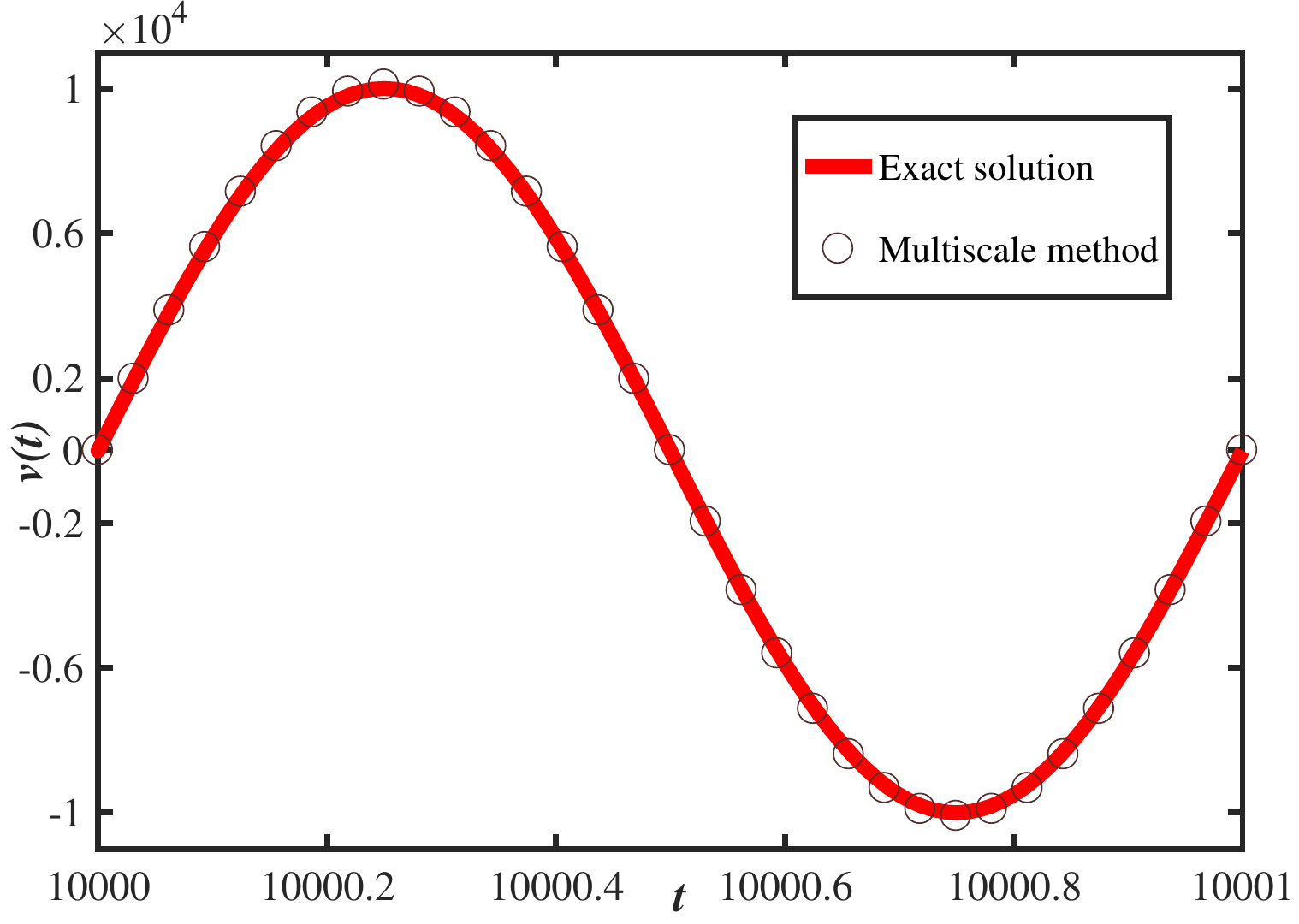}
	\caption{Comparison of the numerical computational results and exact solution for a long time.}
	\label{compare2}
\end{figure}

\begin{figure}[H]
	\centering
	\includegraphics[scale=0.35]{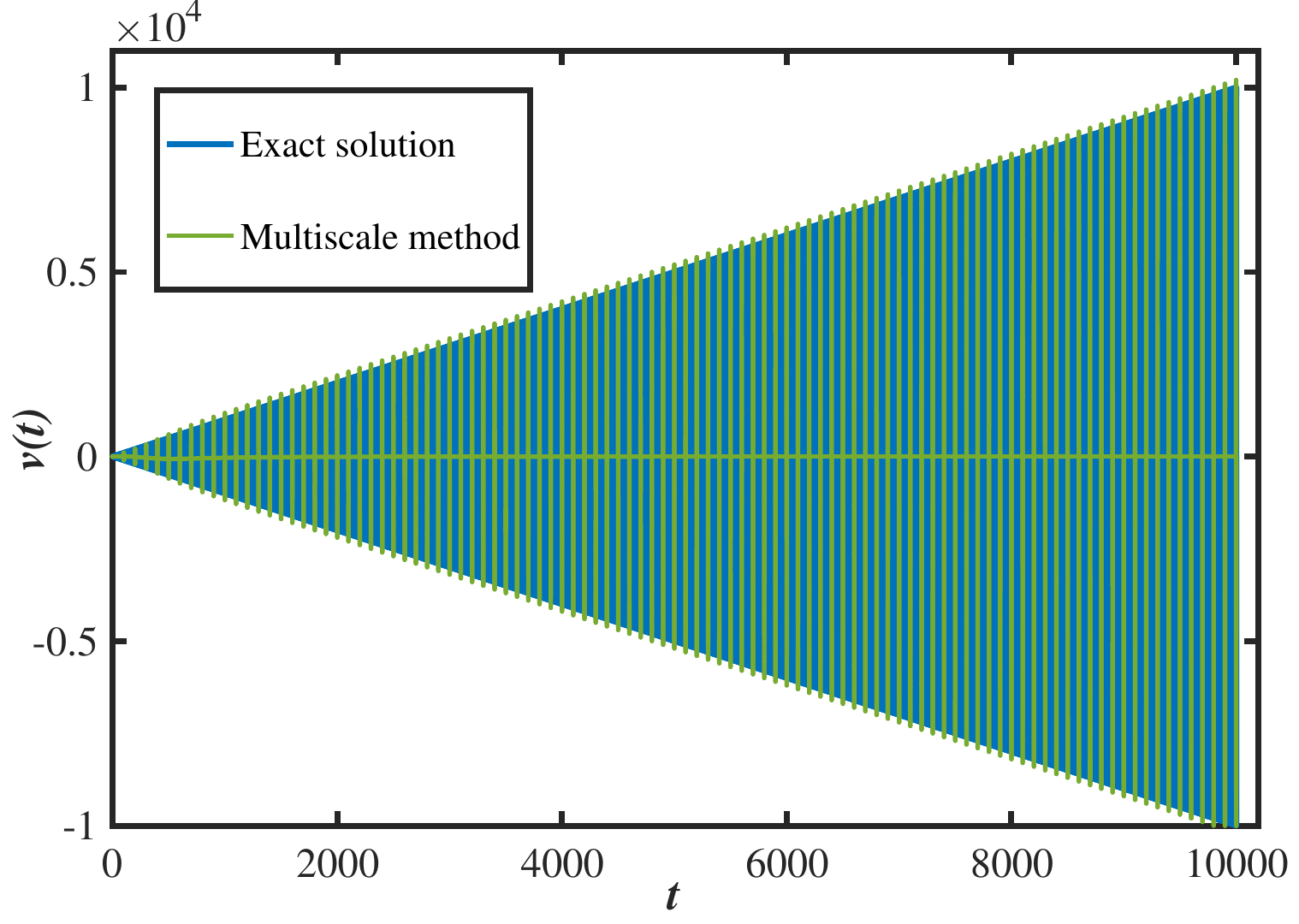}
	\includegraphics[scale=0.35]{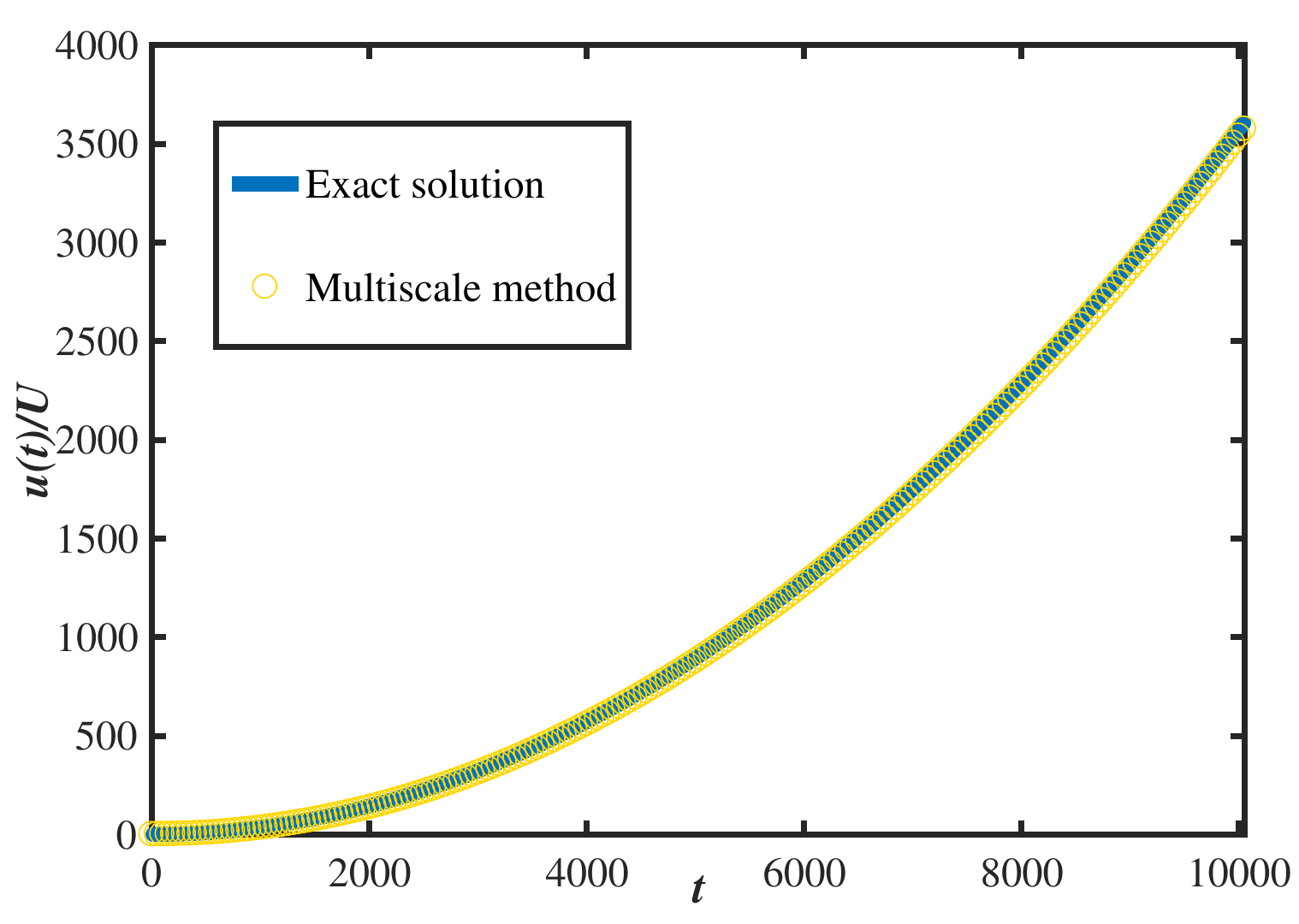}
	\caption{Overall comparison of the numerical results and exact solution of $v(t)$ and $u(t)$.}
	\label{compare3}
\end{figure}

For the discretization of fast variable, we use the first-order implicit Euler scheme. For Algorithm \ref{algorithm11}, we set the time step $\triangle t=1/32$, $\alpha=0.4$, and $\varepsilon=5\times 10^{-5}$ for long-term numerical simulation. For Algorithm \ref{algorithm13}, we set the micro-scale time step $\triangle t=1/32$ and the macro-scale time step $\triangle T=30$ to carry out the simulation. Fig. \ref{compare2} presents the comparison of the numerical results with analytical solution for a long time. Fig. \ref{compare3} shows the calculation results of the multiscale method and the comparison with the exact solution. It is not difficult to see that the multiscale method has good accuracy. 

In order to better test the effectiveness of the multiscale method, we set the micro-scale time step as $\triangle t=1/100$. Quantitative error measurements and computational time are shown in Table \ref{table1}. We observe that the error decreases as the macro-scale time step decreases, and the computational cost increases accordingly. The convergence rates show that our numerical method is first order accurate. It should be noted that due to the memory characteristic of the fractional operator, the macro-scale time step we choose should be much smaller compared to the computation of the integer time multiscale problem. 

Next, we test the impact of the micro-scale time step on the performance of the multiscale method. As shown in Table \ref{table2}, the micro-scale time step has little effect on the performance of the proposed method. The reason may be that the equation for fast variable is of integer order.

For a fully resolved simulation, the $L^\infty$ error is 0.0223, but the time is up to 19 hours. If we need to calculate for a longer time, as shown in the second equation of (\ref{the48}), the calculation of fractional operator will greatly increase the calculation cost, and the multiscale method presents a huge advantage. According to this example, the calculation time is saved by more than 1000 times, and the maximum error does not exceed 4\% ($\triangle T=2$). When the minimum macro-scale time step is $\triangle T=1$, we find that the acceleration is nearly 300 times and the maximum error is only 2\%. It is worth pointing out that 
the discretization scheme of Eq. (\ref{the55}) is of low-order accuracy, which may cause large error results. \\

\begin{table}[H]
	{\footnotesize
		\caption{$L^1$ error norms, $L^\infty$ error norms, convergence rates and CPU time of slow variable $U$ with fixed $\triangle t=1/100$.}\label{table1}
		\begin{center}
			\begin{tabular}{|c|c|c|c|c|c|} \hline  		
				$\triangle T$ & $L^1$ error & order & $L^\infty$ error & order &CPU time \\ \hline
				20 & 7.0964 & \ & 14.2370 & \ & 1.240 s \\  
				10 & 3.5567 & 0.9965 & 7.1271 & 0.9983 & 3.279 s \\  
				5 & 1.7811 & 0.9978 & 3.5666 & 0.9988 & 11.874 s \\  
				2 & 0.7133 & 0.9987 & 1.4276 & 0.9993 & 63.266 s \\
				1 & 0.3568 & 0.9994 & 0.7139 & 0.9998 & 237.580 s \\ \hline
			\end{tabular}
		\end{center}
	}
\end{table}
\begin{table}[H]
	{\footnotesize
		\caption{$L^1$ error norms, $L^\infty$ error norms and CPU time of slow variable $U$ with fixed $\triangle T=1$.}\label{table2}
		\begin{center}
			\begin{tabular}{|c|c|c|c|} \hline  		
				$\triangle t$ & $L^1$ error & $L^\infty$ error &CPU time \\ \hline
				1/16 & 0.3568 & 0.7141 & 238.889 s \\  
				1/32 & 0.3568 & 0.7140 & 231.791 s \\  
				1/64 & 0.3568 & 0.7140 & 255.911 s \\  
				1/128 & 0.3568 & 0.7139 & 261.107 s \\ \hline
			\end{tabular}
		\end{center}
	}
\end{table}	

\textbf{Example 3.} Next we test two Riccati equations coupled with each other
\begin{equation}\label{the58}
\begin{split}
&\frac{d}{dt}v(t)+u(t)v(t)^2+u(t)v(t)=f_3(t)\\
&\frac{d^\alpha}{dt^\alpha}u(t)=\varepsilon\left\{-v(t)u(t)^2+(tsin^2(\pi t)+1)(\varepsilon\Gamma(2-\alpha)t+1)^2+t^{1-\alpha}\right\}\\
&v(0)=1,\quad u(0)=1,\quad t\in [0,8001],
\end{split}
\end{equation}
where $f(t)=sin^2(\pi t)+\pi tsin(2\pi t)+(\varepsilon\Gamma(2-\alpha)t+1)(tsin^2(\pi t)+1)^2+(\varepsilon\Gamma(2-\alpha)t+1)(tsin^2(\pi t)+1)$. We set $\alpha=0.8$ and the scale separation parameter $\varepsilon=5\times10^{-5}$.

The approximate equations are as follow
\begin{equation}\label{the59}
\begin{split}
&\frac{d}{dt}v(t)+Uv(t)^2+Uv(t)=f_3(t)\\
&\frac{d^\alpha}{dt^\alpha}U(t)=\varepsilon\int_t^{t+1}\left\{-v_U(s)U(t)^2+(tsin^2(\pi t)+1)(\varepsilon\Gamma(2-\alpha)t+1)^2+t^{1-\alpha}\right\}ds\\
&v_U([t])=v_U([t]+1), \ U(0)=u_0=1,\quad t\in [0,8001].
\end{split}
\end{equation}
We has exect solutions
\begin{equation}\label{the60}
\begin{split}
v(t)=tsin^2(\pi t)+1,\quad u(t)=\varepsilon\Gamma(2-\alpha)t+1.
\end{split}
\end{equation}

\begin{figure}[H]
	\centering
	\includegraphics[scale=0.35]{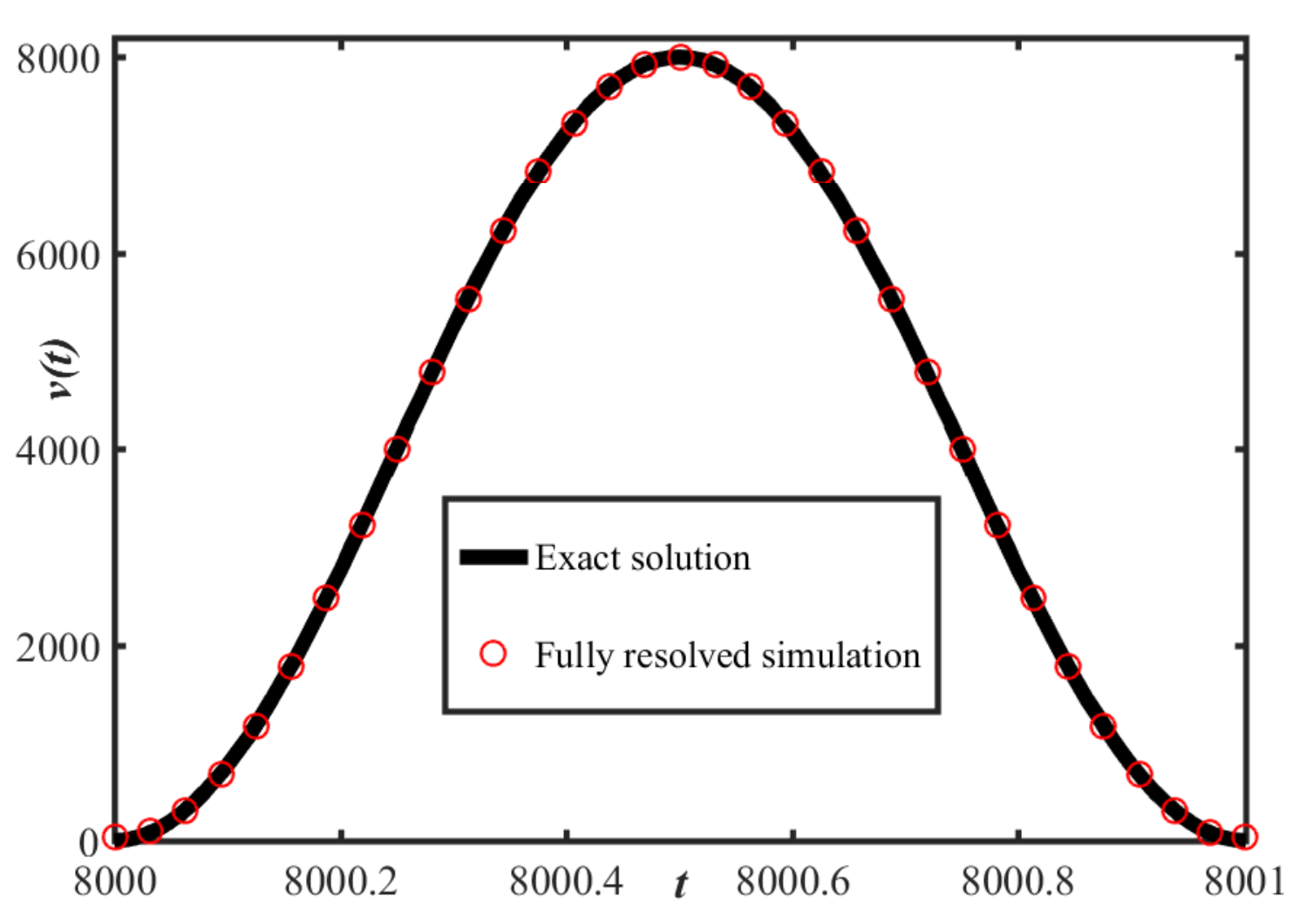}
	\includegraphics[scale=0.35]{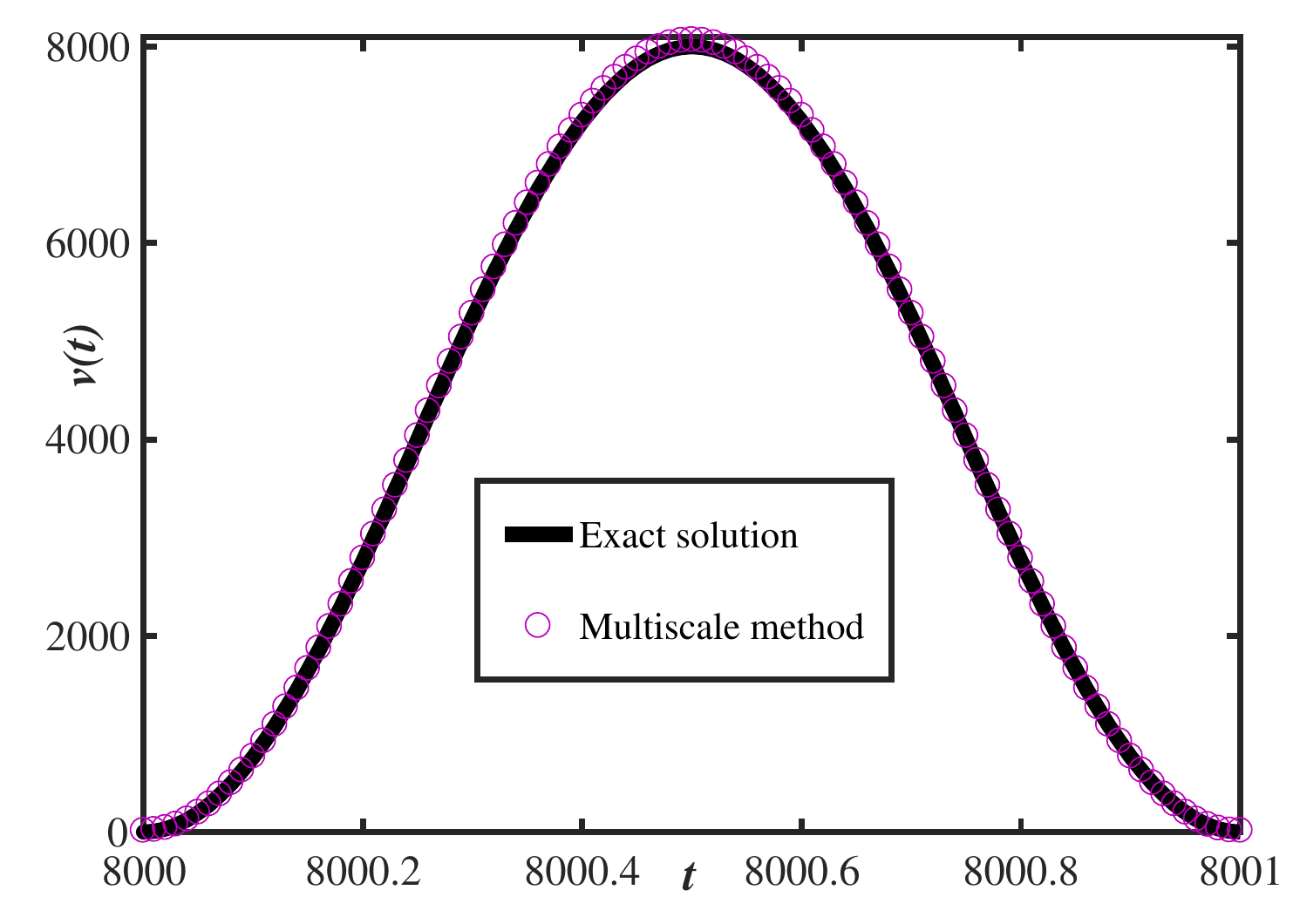}
	\caption{Comparison of the numerical computational results and exact solution for a long time.}
	\label{compare4}
\end{figure}

\begin{figure}[H]
	\centering
	\includegraphics[scale=0.38]{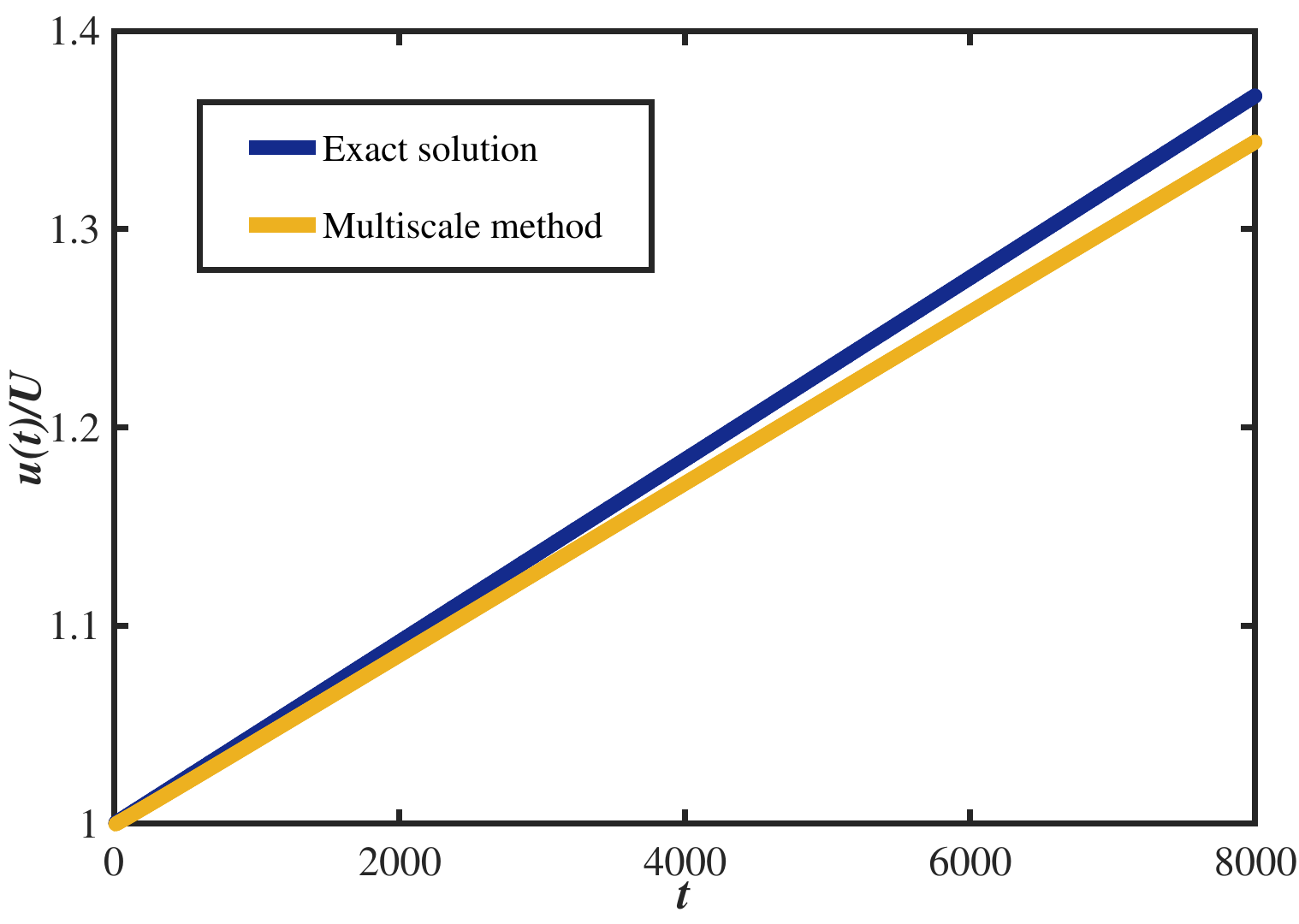}
	\caption{Comparison of the multiscale method and exact solution of $u(t)$ with $\triangle T=10$.}
	\label{compare5}
\end{figure}

We use the frist-order implicit Euler scheme to discretize the fast variable. For the fully resolved simulation, we set the time step $\triangle t=1/32$. For multiscale method, we set the micro-scale time step $\triangle t=1/100$ and the macro-scale time step $\triangle T=10$ for simulation. The comparison between the numerical results and the exact solutions is shown in Fig. \ref{compare4} and \ref{compare5}. We can see that the effect of the multiscale method is good for fast variable and slightly worse for slow variable, which reminds us to decrease the macro-scale time step or apply the high-order discrete scheme. The comparison between the multiscale method and the fully resolved simulation is shown in Table \ref{table3}. The multiscale method also shows superior effects for nonlinear equations, the calculation time is saved by 10000 tiems and the relative error is only 2\%.

Compared with Example 2, we have an interesting conclusion: when the time scale separation becomes larger, the efficiency of the multiscale method is improved since we can select a larger macro-scale time step for computation.

\begin{table}[H]
	{\footnotesize
		\caption{$L^1$ error norms, $L^\infty$ error norms, and CPU time of slow variable $u/U$.}\label{table3}
		\begin{center}
			\begin{tabular}{|c|c|c|c|c|} \hline  		
				$\triangle T$ & $\triangle t$ &$L^1$ error & $L^\infty$ error & CPU time \\ \hline
				Fully resolved simulation & 1/32 & 1.61e-3 & 5.50e-3 & 13.2h \\  
				10 & 1/100 & 1.20e-2 & 2.33e-2 & 3.590 s \\  
				5 & 1/100 & 1.18e-2 & 2.31e-2  & 7.233 s \\  
				2 & 1/100 & 1.17e-2 & 2.30e-2  & 39.366 s \\
				1 & 1/100 & 1.17e-2 & 2.30e-2 & 155.915 s \\ \hline
			\end{tabular}
		\end{center}
	}
\end{table}

\subsection{Approximation of the periodic problem}
An example without analytical solution will be given below to test the effect of the multiscale method.

\ \textbf{Example 4.} In our last example, we solve problem \ref{the1} and set the unknown function as follows
\begin{equation}\label{the61}
\begin{split}
&\frac{d}{dt}v(t)+u(t)^2v(t)^2=t^{1/4}sin(2\pi t)+5\\
&\frac{d^\alpha}{dt^\alpha}u(t)=\varepsilon \left(v(t)u(t)^2 \right)\\
&v(0)=1,\quad u(0)=0.5,\quad t\in [0,10000],
\end{split}
\end{equation}
where $\alpha=0.6$ and $\varepsilon=5\times10^{-5}$.

The approximate periodic equations are as follow
\begin{equation}\label{the62}
\begin{split}
&\frac{d}{dt}v_U(t)+U^2v_U(t)^2=t^{1/4}sin(2\pi t)+5\\
&\frac{d^\alpha}{dt^\alpha}U(t)=\varepsilon U(t)^2\int_t^{t+1}v_{U(t)}(s)ds\\
&v_U([t])=v_U([t]+1), \ U(0)=u_0=0.5,\quad t\in [0,10000].
\end{split}
\end{equation}

\begin{figure}[H]
	\centering
	\includegraphics[scale=0.35]{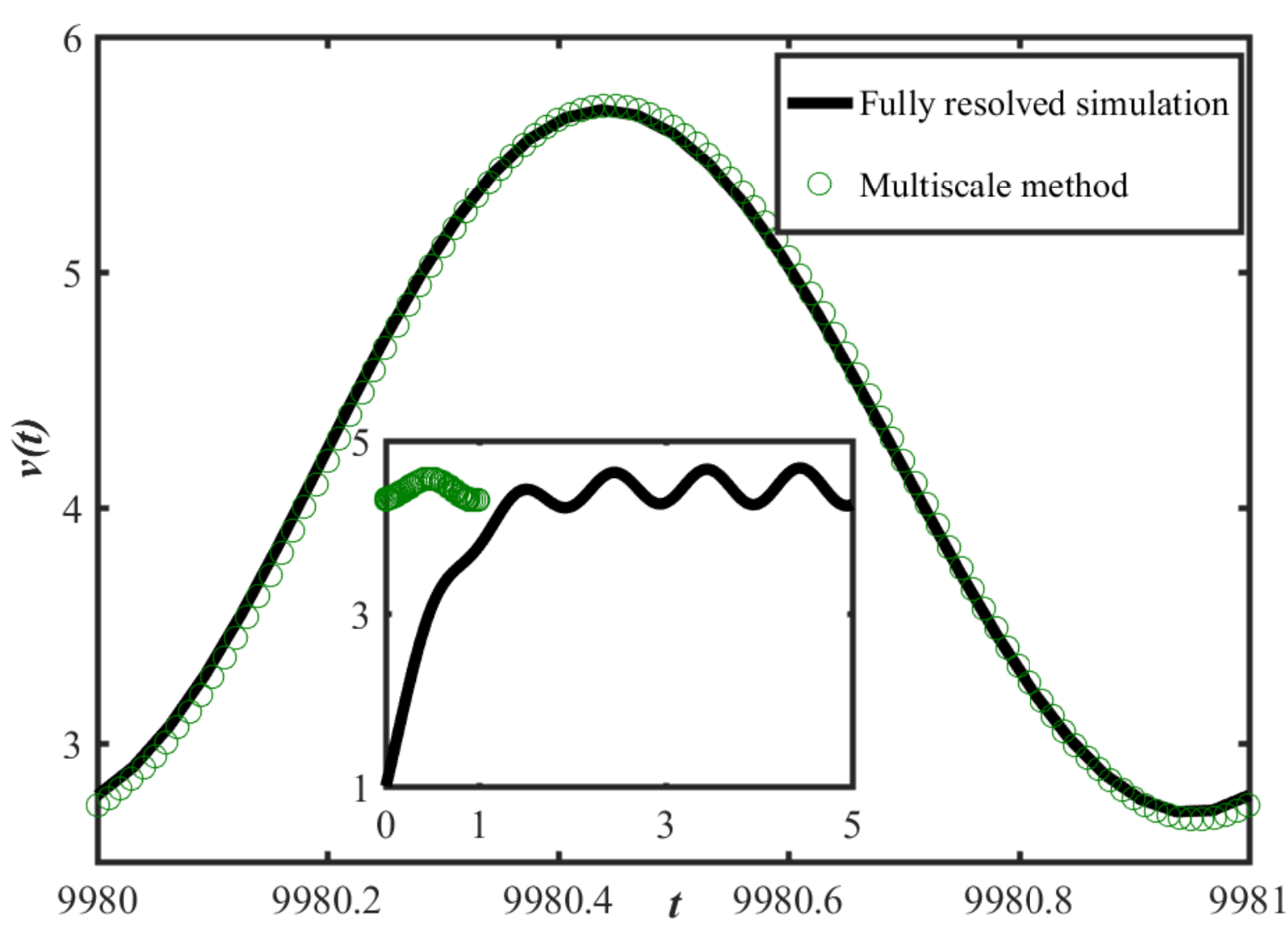}
	\includegraphics[scale=0.35]{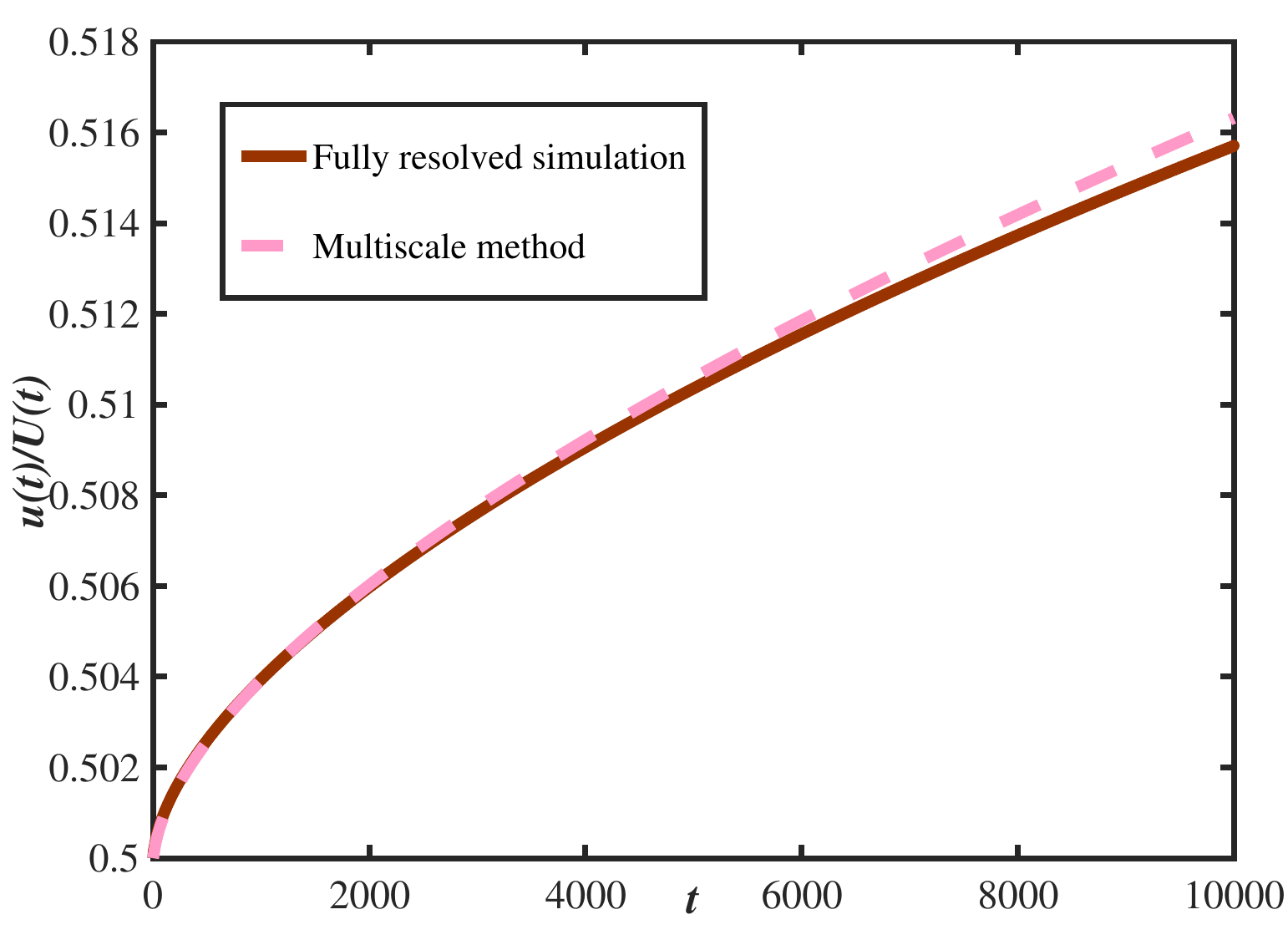}
	\caption{Comparison of the numerical multiscale simulation and the fully resolved  simulation of $v(t)$ and $u(t)$.}
	\label{compare6}
\end{figure}

For multiscale simulation, we set $\triangle t=1/100$ and $\triangle T=20$. For comparison we show the results obtained with the fully resolved simulation for time step $\triangle t=1/32$, the fast variable equation is discretized using the first-order implicit Euler scheme. It can be seen from Fig. \ref{compare6} that the results obtained by applying the multiscale method are close to the fully resolved simulation. It is worth mentioning that the time taken for the fully resolved simulation is 19.8 hours. To compared with the fully resolved simulation, we take $\triangle T=100$, $50$, $10$, $5$ to test the accuracy and computational time of the multiscale method. As shown in Table \ref{table4}, we find that the results obtained by the multiscale method and the direct method have an error limit, which can not be approximated by decreasing the macro/micro-scale time step. This result is not difficult to understand. In addition to the error analyzed in Section \ref{section3}, the error limit should also include a numerical error.

\begin{table}[H]
	{\footnotesize
		\caption{$L^1$ error norms, $L^\infty$ error norms and CPU time of slow variable $u/U$ with fixed $\triangle t=1/100$.}\label{table4}
		\begin{center}
			\begin{tabular}{|c|c|c|c|} \hline  		
				$\triangle T$ & $L^1$ error & $L^\infty$ error &CPU time \\ \hline
				100 & 1.893e-4 & 5.720e-4 & 0.235 s \\  
				50 & 2.275e-4 & 6.012e-4 & 0.422 s \\   
				10 & 2.598e-4 & 6.249e-4 & 6.838 s \\  
				5 & 2.640e-4 & 6.279e-4 & 11.761 s \\ \hline
			\end{tabular}
		\end{center}
	}
\end{table}

\section{Conclusions and remarks}
\label{section6}
In this paper, a nonlinear system of coupled fractional ordinary differential equations with a periodic applied force and a temporal multiscale feature is studied. In order to improve the computational efficiency, we derive an auxilary local periodic equation so as to formulate a multiscale approximate problem, and prove that the error of its solutions to the original solution is $O(\varepsilon)$. The finite difference method is applied to solve the original and approximate multiscale equations. Several numerical examples show the superior performance of the multiscale method for the calculation of fractional differential equations. Our analytical process and numerical method can be used as a general framework for the simulation of more practical model equations with temporal multiscale feature. It should be noted that currently this multiscale method is theoretically established  only for problems with an applied periodic force in the fast variable equation (e.g. mimicking the periodic pulse in a blood flow). 

We focus on the effect of the change of macro-scale time step on the efficiency of multiscale method, as the micro-scale time step has little effect on it. When the time scales of fast and slow variables are more separated, we can use a larger macro-scale time step without much influence of the accuracy of the multiscale method, therefore, computational time can be largely saved.

\section*{Acknowledgments}
Z. Wang would like to thank Dr. Weifeng Zhao for discussions during this work. This work is supported by the National Natural Science Foundation of China (Nos. 11861131004 and 11771040) and the Fundamental Research Funds for the Central Universities (No. 06500073).

\bibliography{Ref}

\end{document}